\documentclass[11pt]{article}
\usepackage{amsmath,amssymb,amsthm}
\usepackage[a4paper,margin=25mm]{geometry}
\usepackage{microtype}
\usepackage[hidelinks]{hyperref}
\newtheorem{theorem}{Theorem}
\newtheorem{lemma}[theorem]{Lemma}

\theoremstyle{definition}
\newtheorem{problem}{Problem}
\theoremstyle{plain}
\newcommand{\PP}{\mathcal P}
\newcommand{\A}{\mathcal A_3}
\newcommand{\B}{\mathcal B_3}
\newcommand{\C}{\mathcal C_3}
\newcommand{\LL}{\mathcal L}
\DeclareMathOperator{\vol}{vol}
\title{Density and distribution of quotients associated with Goldring's $w$-function}

\author{%
	Quan-Hui Yang\textsuperscript{a},\quad
	Jun-Jia Zhao\textsuperscript{b,}%
	\thanks{Corresponding author.\protect\\
		E-mail addresses:
		\href{mailto:yangquanhui01@163.com}
		{yangquanhui01@163.com} (Q.-H. Yang),
		\href{mailto:jjzhao1998@163.com}
		{jjzhao1998@163.com} (J.-J. Zhao).}%
	\\[0.8em]
	\parbox{0.91\textwidth}{%
		\small\itshape\raggedright
		\textsuperscript{a} Ministry of Education
		Key Laboratory for NSLSCS and School of
		Mathematical Sciences, Nanjing Normal
		University, Nanjing 210023, China
		\\[0.35em]
		\textsuperscript{b} School of Mathematical
		Sciences, Nanjing Normal University of
		Special Education, Nanjing 210038, PR China%
	}%
}
\date{}

\begin{document}
	\maketitle
	\begin{abstract}
		For primes $p,q,r$ that are not all equal, let
		$n=pqr$ and define Goldring's $w$-function by
		$w(n)=P(p+q)P(p+r)P(q+r)$, where $P(a)$ denotes
		the largest prime factor of $a$. In this paper,
		we resolve three problems of Chen, Shi and Wu.
		As $n\to\infty$ through such products, the sets
		of finite limit points of $P(w(n))/P(n)$ and
		$w(n)/n$ are $[0,1]$ and $[0,\infty)$,
		respectively.
		The same conclusions hold when $n$ is restricted
		either to $p^2q$ with distinct primes $p,q$ or
		to squarefree products of three primes.
		Nevertheless, on each of these three sets,
		the proportion of integers $n\le X$ satisfying
		$P(w(n))/P(n)>\alpha$ tends to zero as
		$X\to\infty$, for every fixed $\alpha>0$.
		We also determine the complete set of finite
		joint limit points of the two quotients for
		$n=p^2q$ with $p\ne q$.
	\end{abstract}

	

	\section{Introduction}
	Let $\PP$ denote the set of primes. For a positive integer $a$,
	let $P(a)$ be its largest prime factor, with $P(1)=1$. We consider
	\[
	\A=\{pqr:p,q,r\in\PP\}\setminus\{p^3:p\in\PP\},\qquad
	\B=\{p^2q:p,q\in\PP,\ p\ne q\},\qquad \C=\A\setminus\B.
	\]
	Thus $\A=\B\sqcup\C$: the integers in $\B$ have one repeated
	prime factor, while those in $\C$ are squarefree products of three primes.
	
	Goldring~\cite{Goldring} introduced the map
	\[
	w(pqr)=P(p+q)P(p+r)P(q+r)
	\]
	and proved $w(\A)\subset\A$. We study the quotients
	\[
	R(n)=\frac{P(w(n))}{P(n)},\qquad S(n)=\frac{w(n)}n.
	\]
	The first compares the largest prime factors of $w(n)$ and $n$,
	while the second compares their sizes. The two ratios can behave differently.
	For example, $63=3^2\cdot7$ and
	\[
	w(63)=P(6)P(10)^2=75,\qquad
	R(63)=\frac57<1,\qquad S(63)=\frac{25}{21}>1.
	\]
	
	\subsection{Questions and main results}
	Chen, Shi and Wu~\cite[Corollary 5]{CSW} established the following extremal limits:
	\[
	\begin{aligned}
		\liminf_{\substack{n\to\infty\\n\in\A}}R(n)&=0,&
		\limsup_{\substack{n\to\infty\\n\in\A}}R(n)&=1,\\
		\liminf_{\substack{n\to\infty\\n\in\A}}S(n)&=0,&
		\limsup_{\substack{n\to\infty\\n\in\A}}S(n)&=\infty.
	\end{aligned}
	\]
	We restate Problems~7--9 of Chen, Shi and Wu~\cite{CSW}
	as Problems~1--3 below. The first two ask whether the sets of
	quotient values are dense in the indicated intervals.

	\begin{problem}
		Is $R(\A)$ dense in $[0,1]$?
	\end{problem}
	
	\begin{problem}
		Is $S(\A)$ dense in $[0,\infty)$?
	\end{problem}
	
	The third question concerns relative frequency
	rather than the existence of approximating values.
	For a set
	$\mathcal S$ of positive integers, write
	$\mathcal S(X)=|\mathcal S\cap[1,X]|$. For $0<\alpha\le1$, define
	\[
	D_{\A}(\alpha)=\lim_{X\to\infty}
	\frac{|\{n\le X:n\in\A,\ R(n)\le\alpha\}|}{\A(X)},
	\]
	when the limit exists. Chen, Shi and Wu~\cite[Theorem 6]{CSW} proved that
	$D_{\A}(\alpha)=1$ for $1/2<\alpha\le1$, and that
	\[
	|\{n\le X:n\in\A,\ R(n)\le\alpha\}|\gg\alpha^2\A(X)
	\qquad(0<\alpha\le1/2),
	\]
	with an absolute implied constant. They further posed the following question.
	
	\begin{problem}
		Does $D_{\A}(\alpha)$ exist for every
		$0<\alpha\le1/2$? If so, what is its value?
	\end{problem}
	
	We answer these three questions, and also consider the corresponding
	assertions on $\B$ and $\C$. Throughout the paper, a limit point means
	a limit along a sequence of integers $n\to\infty$. We use the standard
	notation $f\gg g$ for a lower bound $f\ge c g$ with a positive constant
	$c$, and $f\ll g$ for the corresponding upper bound. A subscript
	indicates permitted dependence of the constant. Also, $f\sim g$ means
	$f/g\to1$, while $f=o(g)$ means $f/g\to0$.
	
	Our first result determines the individual finite limit-point sets and gives
	quantitative lower bounds for the number of integers whose selected quotient
	lies in a prescribed open interval.
	\begin{theorem}\label{th:dense}
		For each $\mathcal S\in\{\B,\C\}$, the sets of finite limit points of
		$R(n)$ and $S(n)$, with $n\in\mathcal S$, are $[0,1]$ and
		$[0,\infty)$, respectively.
		More precisely, if $I\subset(0,1)$ and $J\subset(0,\infty)$ are
		nonempty open intervals, then, for all sufficiently large $X$,
		\begin{align*}
			|\{n\le X:n\in\B,\ R(n)\in I\}|&\gg_I X^{2/3}(\log X)^{-3},\\
			|\{n\le X:n\in\B,\ S(n)\in J\}|&\gg_J X^{2/3}(\log X)^{-3},\\
			|\{n\le X:n,w(n)\in\C,\ R(n)\in I\}|&\gg_I X(\log X)^{-6},\\
			|\{n\le X:n,w(n)\in\C,\ S(n)\in J\}|&\gg_J X(\log X)^{-6}.
		\end{align*}
		The constants and the thresholds may depend on the indicated interval.
	\end{theorem}
	
	Each lower bound concerns one quotient at a time; the bounds for
	$\C$ also require $w(n)\in\C$. Since these bounds tend to infinity, every such
	interval contains values of the corresponding quotient for arbitrarily large $n$.
	
	The answer to the frequency question is that $R(n)$ tends to zero in
	relative density, even on each of the two parts of $\A$.
	\begin{theorem}\label{th:distribution}
		For every fixed $\alpha>0$ and every
		$\mathcal S\in\{\A,\B,\C\}$, we have
		\begin{equation}\label{eq:density}
			|\{n\le X:n\in\mathcal S,\ R(n)>\alpha\}|=o(\mathcal S(X)).
		\end{equation}
		In particular, $D_{\A}(\alpha)=1$ for every $0<\alpha\le1$.
	\end{theorem}
	
	For $\B$, we use the fixed-shift sieve estimate of Chen and
	Shi~\cite{CSparents},
	truncating the sum over the repeated prime. We include the deduction
	for completeness. For $\C$, we use an estimate uniform in the prime
	shift; the result for $\A$ then follows by combining the two cases.
	
	Thus every value in $[0,1]$ can be approximated along a sequence of
	integers $n\to\infty$, although $R(n)$ converges to zero in relative density.
	However, $R(n)$ does not tend to zero along every sequence
	$n\to\infty$ in $\B$ or $\C$: for $0<a<b<1$ and
	$\mathcal S\in\{\B,\C\}$, the two theorems give simultaneously
	\[
	\#\{n\in\mathcal S:a<R(n)<b\}=\infty,\qquad
	\frac{\#\{n\le X:n\in\mathcal S,\ a<R(n)<b\}}{\mathcal S(X)}
	\longrightarrow0.
	\]
	
	Knowledge of the individual limit-point sets of the two quotients
	does not determine their joint limit-point set. On $\B$, the two quotients are linked by a simple
	identity already used by Goldring~\cite[formula (14)]{Goldring}.
	The three pairwise sums for $p,p,q$ are $2p,p+q,p+q$.
	Since $P(2p)=p$ for every prime $p$, including $p=2$, we have
	\[
	w(p^2q)=P(2p)P(p+q)^2=pP(p+q)^2.
	\]
	For distinct odd primes $p,q$, put $s=P(p+q)$. We have
	$2s\mid p+q$: for $s>2$, this follows from
	$2\mid p+q$, $s\mid p+q$ and $(2,s)=1$; for $s=2$,
	$p+q=2^h\ge8$, so $4\mid p+q$. Hence
	\[
	u=\frac pq>0,\qquad
	m=\frac{p+q}{2s}\in\mathbb Z_{\ge1},\qquad
	\frac sq=\frac{u+1}{2m}.
	\]
	It follows that
	\[
	R(p^2q)=\frac{\max\{p,s\}}{\max\{p,q\}},\qquad
	S(p^2q)=\frac{s^2}{pq}.
	\]
	For $m\ge1$ and $u>0$, define the coordinate functions and their
	parametrized curves by
	\begin{equation}\label{eq:curves}
		\begin{aligned}
			r_m(u)&=\frac{\max\bigl\{u,\frac{u+1}{2m}\bigr\}}{\max\{u,1\}},
			&g_m(u)&=\frac{(u+1)^2}{4m^2u},\\[3pt]
			F_m(u)&=\bigl(r_m(u),g_m(u)\bigr),
			&\Gamma_m&=\{F_m(u):u>0\}.
		\end{aligned}
	\end{equation}
	Thus $(R(p^2q),S(p^2q))=F_m(u)$ for the above integer $m$ and
	positive real $u$.
	
	\begin{theorem}\label{th:joint}
		The set $\LL_{\B}$ of all finite limits of $(R(n),S(n))$ along
		$n\in\B$, $n\to\infty$, is exactly
		\begin{equation}\label{eq:joint}
			\LL_{\B}=([0,1]\times\{0\})\cup(\{0\}\times[0,\infty))
			\cup\bigcup_{m\ge1}\Gamma_m.
		\end{equation}
		The same set is obtained if both distinct prime divisors of $n$
		are required to be odd.
	\end{theorem}
	
	For each fixed $m$, the parametrization yields the curve $\Gamma_m$.
	The axis components arise as limits along suitable sequences with $m\to\infty$.
	For example, the section at first coordinate $1/2$ is
	\[
	\LL_{\B}\cap(\{1/2\}\times[0,\infty))
	=\{(1/2,0)\}\cup
	\left\{\left(\frac12,\frac9{8m^2}\right):m\ge2\right\}.
	\]
	Thus $(1/2,1)$ is not a joint limit, although $1/2$ is a limit of $R$
	and $1$ is a limit of $S$ separately. The joint limit set is therefore
	a proper subset of $[0,1]\times[0,\infty)$.
	Theorems~\ref{th:dense} and~\ref{th:distribution} answer Problems 7--9
	of~\cite{CSW}.
	
	\subsection{Earlier work and the methods used here}
	Previous studies of Goldring's $w$-function have emphasized iteration,
	inverse images, and comparisons between an input and its image.
	Goldring~\cite{Goldring} proved that every orbit in $\A$ eventually
	enters the cycle
	\[
	20\longrightarrow98\longrightarrow63\longrightarrow75
	\longrightarrow20.
	\]
	Writing $\operatorname{ind}(n)$ for the number of steps needed to
	reach this cycle, Chen and Shi~\cite{CSindex} improved Goldring's
	bound $\operatorname{ind}(n)\ll\pi(P(n))$ to
	$\operatorname{ind}(n)\ll(\log P(n))^2$ and, using the
	Green--Tao theorem on prime progressions~\cite{GTAP}, proved that the
	indices are unbounded. Extensions to functions on products of more
	prime factors were studied by Shi~\cite{Shi} and Shi and
	Yang~\cite{SY}.
	
	For the inverse problem, a parent of $n$ is an integer $m$ with
	$w(m)=n$. Chen and Shi~\cite{CSparents} proved that, for every fixed
	$h\ge1$, infinitely many members of $\B$ have at least $h$ parents
	in $\B$, while infinitely many have none. Quantitative refinements for $\B$ and $\C$ were obtained
	by Jia~\cite{Jia} and Chen, Shi and Wu~\cite{CSW}. For size comparisons, Chen, Shi and Wu~\cite[Theorem 4]{CSW}, using
	estimates of Baker and Harman~\cite{BH}, established
	\[
	\begin{aligned}
		\liminf_{\substack{n\to\infty\\n\in\A}}
		\frac{\log P(w(n))}{\log P(n)}&\le0.2962, &
		\liminf_{\substack{n\to\infty\\n\in\A}}
		\frac{\log w(n)}{\log n}&\le0.5924,
	\end{aligned}
	\]
	and
	\[\limsup_{\substack{n\to\infty\\n\in\A}}
	\frac{\log w(n)}{\log n}\ge1.354.\]
	These estimates concern the power scale: for example,
	\[
	\frac{\log P(w(n))}{\log P(n)}\le\theta
	\quad\Longleftrightarrow\quad P(w(n))\le P(n)^\theta.
	\]
	They do not prescribe positive finite values of $R(n)$ or $S(n)$.
	More recently, Lichtman~\cite[Theorem 1.1]{Lichtman} proved that for
	each fixed nonzero integer $a$ and every $\beta>15/(32\sqrt e)$,
	infinitely many primes $p$ satisfy $P(p-a)\le p^\beta$.
	With $a=-2$ and $w(4p)=2P(p+2)^2$ for odd $p$, this implies
	\[
	\frac{\log P(w(4p))}{\log P(4p)}\le\beta,
	\qquad \frac{\log w(4p)}{\log(4p)}\le2\beta+o(1).
	\]
	Letting $\beta\downarrow15/(32\sqrt e)$ improves the first two
	numerical upper bounds above to $15/(32\sqrt e)$ and
	$15/(16\sqrt e)$, respectively.
	
	For the existence results, we apply the unconditional
	complexity-at-most-two case of Green and Tao's theorem on linear
	equations in primes~\cite{GT}, verifying local solvability for
	explicit systems of three and six forms. For relative density, we use
	the fixed-shift estimate of Chen and Shi~\cite{CSparents}
	for $\B$, a uniform prime-shift estimate for $\C$, and the counting
	formulas of Chen, Shi and Wu~\cite[Lemma 6.1]{CSW}.
	Standard background on prime distribution and sieve methods may be
	found in~\cite[Chapters 2, 3 and 11]{MV}.
	
	Section~2 gives the prime-pattern and sieve estimates. Section~3
	proves the separate limit-point, joint limit-point and relative-density
	results.
	
	\section{Preliminaries}
	Write $L_i(\boldsymbol x)=v_i\cdot\boldsymbol x$, where
	$v_i\in\mathbb Z^d\setminus\{0\}$. For homogeneous forms, complexity
	at most $s$ means that, for every $i$, there is a partition
	\[
	\{1,\ldots,t\}\setminus\{i\}
	=\coprod_{\nu=0}^{s}I_{i,\nu},\qquad
	v_i\notin\operatorname{span}_{\mathbb Q}
	\{v_j:j\in I_{i,\nu}\}\quad(0\le\nu\le s);
	\]
	see~\cite[Definition 1.5]{GT}. We shall use the criteria
	\[
	\begin{aligned}
		v_i\notin\mathbb Qv_j
		&\quad\Longleftrightarrow\quad v_i,v_j
		\text{ are not proportional},\\
		\det(v_i,v_j,v_k)\ne0
		&\quad\Longrightarrow\quad
		v_i\notin\operatorname{span}_{\mathbb Q}\{v_j,v_k\}
		\qquad(d=3).
	\end{aligned}
	\]
	
	We use Green and Tao~\cite[Corollary 1.7 and (1.8)]{GT} as an
	external theorem. Let the nonzero integral homogeneous forms
	$L_1,\ldots,L_t$ and a compact convex body $K\subset\mathbb R^d$
	be fixed. Assume
	\[
	\begin{gathered}
		\operatorname{complexity}(L_1,\ldots,L_t)\le2,\qquad
		\operatorname{int}(K)\ne\varnothing,\qquad
		L_i(\boldsymbol x)>0\quad(\boldsymbol x\in K,\ 1\le i\le t),\\
		\text{for every prime }\ell,\quad
		\exists\boldsymbol a\in(\mathbb Z/\ell\mathbb Z)^d:
		\prod_{i=1}^tL_i(\boldsymbol a)\not\equiv0\pmod\ell.
	\end{gathered}
	\]
	The last condition is local solvability. Then, with
	$NK=\{N\boldsymbol x:\boldsymbol x\in K\}$,
	\begin{equation}\label{eq:GT}
		|\{\boldsymbol a\in NK\cap\mathbb Z^d:
		L_i(\boldsymbol a)\in\PP\ (1\le i\le t)\}|
		\sim c_{L,K}\frac{N^d}{(\log N)^t},\qquad c_{L,K}>0.
	\end{equation}
	This complexity-at-most-two case is unconditional. The constants are
	\[
	\begin{aligned}
		c_{L,K}&=\vol(K)\prod_\ell\beta_\ell,\\
		\beta_\ell&=\left(1-\frac1\ell\right)^{-t}
		\frac{|\{\boldsymbol a\bmod\ell:
			\ell\nmid\prod_iL_i(\boldsymbol a)\}|}{\ell^d}.
	\end{aligned}
	\]
	We briefly justify the positivity of the singular series and the
	normalization in~\eqref{eq:GT}. Local solvability ensures
	$\beta_\ell>0$. Finite complexity implies that no two coefficient
	vectors are proportional. For all but finitely many primes $\ell$,
	the congruence hyperplanes
	$H_i=\{\boldsymbol a\bmod\ell:L_i(\boldsymbol a)=0\}$ therefore satisfy
	\[
	|H_i|=\ell^{d-1},\qquad
	|H_i\cap H_j|=\ell^{d-2}\quad(i\ne j).
	\]
	Inclusion--exclusion (or the corresponding single-hyperplane count
	when $t=1$) gives
	\[
	\frac{\bigl|\{\boldsymbol a\bmod\ell:
		\ell\nmid\prod_iL_i(\boldsymbol a)\}\bigr|}{\ell^d}
	=1-\frac{t}{\ell}+O_L(\ell^{-2}),
	\qquad \beta_\ell=1+O_L(\ell^{-2}).
	\]
	Consequently the positive Euler factors have a convergent product,
	$0<\prod_\ell\beta_\ell<\infty$; see also~\cite[Lemma 1.3]{GT}.
	Since $K$ has nonempty interior, its volume is positive.
	
	To match the formulation of~\cite[(1.8)]{GT}, choose fixed $B\ge1$
	with $K\subset[-B,B]^d$ and put $M=\lceil BN\rceil\asymp N$.
	The region and error conditions then hold at this scale. By
	homogeneity, compactness and positivity, uniformly for
	$\boldsymbol x\in NK$ and $1\le i\le t$,
	\[
	cN\le L_i(\boldsymbol x)\le CN,\qquad
	\log L_i(\boldsymbol x)=\log N+O_{L,K}(1)
	\]
	with fixed $c,C>0$. Changing variables $\boldsymbol x=N\boldsymbol y$
	gives the integral main term
	\[
	\int_{NK}\frac{d\boldsymbol x}{\prod_{i=1}^t\log L_i(\boldsymbol x)}
	=\left(\vol(K)+O_{L,K}\!\left(\frac1{\log N}\right)\right)
	\frac{N^d}{(\log N)^t}.
	\]
	This proves the normalization in~\eqref{eq:GT}.
	In particular, the centre of a fixed real box need not have
	integer coordinates.
	
	We will also discard finitely many hyperplanes corresponding to
	equalities between distinct selected forms. For each proper affine
	hyperplane $H$ with fixed coefficients,
	\[
	\#(H\cap NK\cap\mathbb Z^d)
	=O_K(N^{d-1})
	=o\!\left(\frac{N^d}{(\log N)^t}\right).
	\]
	Thus removing finitely many such hyperplanes does not change the
	asymptotic in~\eqref{eq:GT}.
	
	\begin{lemma}\label{lem:three}
		For every positive integer $m$, the three forms
		$a,b,2mb-a$ are locally soluble at every prime and have complexity
		at most one.
	\end{lemma}
	\begin{proof}
		For local solvability modulo $\ell$, take $b=1$.
		If $\ell=2$, choose $a=1$; if $\ell>2$, choose
		$a\not\equiv0,2m\pmod\ell$, which is possible since at most two
		residue classes are excluded. All three forms are then nonzero modulo $\ell$.
		Their coefficient vectors $(1,0)$, $(0,1)$ and $(-1,2m)$ have
		pairwise determinants $1$, $2m$ and $1$, respectively. Thus no two
		vectors are proportional, and partitioning the other two forms into
		singleton classes proves complexity at most one.
	\end{proof}
	
	\begin{lemma}\label{lem:six}
		For every positive integer $m$, the six forms
		\begin{align*}
			L_1&=x,&L_2&=y,&L_3&=z,\\
			L_4&=mx+my-z,&L_5&=mx-my+z,&L_6&=-mx+my+z
		\end{align*}
		are locally soluble at every prime and have complexity at most two.
	\end{lemma}
	\begin{proof}
		Take $(x,y,z)=(a,a,1)$ modulo $\ell$, where
		\[
		a\ne0,\qquad 2ma\ne1,\qquad
		(L_1,\ldots,L_6)=(a,a,1,2ma-1,1,1).
		\]
		For $\ell=2$, choose $a=1$. For $\ell>2$,
		\[
		|\{a\in(\mathbb Z/\ell\mathbb Z)^\times:2ma=1\}|\le1,
		\qquad (\ell-1)-1\ge1,
		\]
		so an admissible $a$ exists. All six values are nonzero.
		
		The coefficient vectors are
		\begin{align*}
			v_1&=(1,0,0),&v_2&=(0,1,0),&v_3&=(0,0,1),\\
			v_4&=(m,m,-1),&v_5&=(m,-m,1),&v_6&=(-m,m,1).
		\end{align*}
		For each $i$, the table partitions the other five forms into three
		classes. The last two columns give $\det(v_i,v_j,v_k)$ for the
		ordered pair $j<k$ in the second and third classes, respectively.
		\[
		\begin{array}{c|ccc|rr}
			i & I_{i,0} & I_{i,1} & I_{i,2}
			& D_{i,1} & D_{i,2}\\ \hline
			1&\{L_2\}&\{L_3,L_4\}&\{L_5,L_6\}&-m&-2m\\
			2&\{L_1\}&\{L_3,L_4\}&\{L_5,L_6\}&m&-2m\\
			3&\{L_1\}&\{L_2,L_5\}&\{L_4,L_6\}&-m&2m^2\\
			4&\{L_1\}&\{L_2,L_3\}&\{L_5,L_6\}&m&-4m^2\\
			5&\{L_1\}&\{L_2,L_3\}&\{L_4,L_6\}&m&4m^2\\
			6&\{L_1\}&\{L_2,L_3\}&\{L_4,L_5\}&-m&-4m^2
		\end{array}
		\]
		Here $D_{i,\nu}$ is the determinant with first row $v_i$
		and the other two rows the coefficient vectors of the forms in
		$I_{i,\nu}$, in increasing index order.

		Since $m\ge1$, every listed determinant is nonzero; hence no pair
		spans its selected vector. For the singleton classes,
		\[
		v_1\notin\mathbb Qv_2,\qquad
		((v_i)_2,(v_i)_3)\ne(0,0)
		\ \Longrightarrow\ v_i\notin\mathbb Qv_1\quad(2\le i\le6).
		\]
		Thus every row satisfies all three span conditions, proving complexity
		at most two by~\cite[Definition 1.5]{GT}.
	\end{proof}
	
	We prove the following bound uniformly in the prime shift $a$,
	as required for the squarefree case below.
	\begin{lemma}\label{lem:shift}
		Fix a positive integer $k$. For $a\in\PP$, let
		\[
		E_k(T;a)=|\{r\le T:r\in\PP,\ (r+a)/k\in\PP\}|.
		\]
		Then $\sup_{a\in\PP}E_k(T;a)=o(T/\log T)$ as $T\to\infty$.
	\end{lemma}
	\begin{proof}
		Fix $z\ge3$, sufficiently large that $W>1$, where
		\[
		W=\prod_{\substack{3\le\ell\le z\\\ell\in\PP,\ \ell\nmid k}}\ell,
		\qquad
		\eta_k(z)=\prod_{\substack{3\le\ell\le z\\\ell\in\PP,\ \ell\nmid k}}
		\frac{\ell-2}{\ell-1}.
		\]
		Throughout the prime-counting argument, $k,z,W$ are fixed.
		
		\smallskip\noindent
		\emph{Step 1.}
		Write $s=(r+a)/k$. Since $r=ks-a$ is determined by $s$,
		\[
		|\{r\text{ counted by }E_k(T;a):r\le z\text{ or }s\le z\}|
		\le\pi(z)+\pi(z)=2\pi(z).
		\]
		For the remaining pairs, $r,s>z$. If $\ell\mid W$, then
		\[
		\ell\nmid r,\qquad
		\ell\mid r+a=ks\ \Longrightarrow\ \ell\mid s,
		\]
		where the latter is impossible because $s$ is prime and $s>z\ge\ell$.
		Thus $r\not\equiv0,-a\pmod\ell$. Since $a$ is prime,
		the number of allowed nonzero classes modulo $\ell$ is
		\[
		h_{a,\ell}=
		\begin{cases}
			\ell-2,&\ell\ne a,\\
			\ell-1,&\ell=a.
		\end{cases}
		\]
		Let $\mathcal R_a$ denote the corresponding set of reduced residue classes modulo $W$, and write
		$h_a = |\mathcal R_a|$.
		By the Chinese remainder theorem,
		\[
		h_a=\prod_{\ell\mid W} h_{a,\ell}.
		\]
		Moreover,
		\[
		\frac{h_a}{\varphi(W)}
		=
		\begin{cases}
			\displaystyle \prod_{\ell\mid W}\frac{\ell-2}{\ell-1}
			= \eta_k(z), & a\nmid W,\\[10pt]
			\displaystyle \prod_{\substack{\ell\mid W\\ \ell\ne a}}
			\frac{\ell-2}{\ell-1}
			= \eta_k(z)\frac{a-1}{a-2}, & a\mid W.
		\end{cases}
		\]
		Since $a\mid W$ implies $a\ge 3$, we have
		\[
		\frac{a-1}{a-2}\le 2,
		\]
		and therefore
		\[
		\frac{h_a}{\varphi(W)} \le 2\eta_k(z).
		\]
		
		\smallskip\noindent
		\emph{Step 2.}
		Let $\pi(T;W,b)$ count primes $r\le T$ with $r\equiv b\pmod W$.
		For $(b,W)=1$, the prime number theorem in arithmetic progressions
		for the fixed modulus $W$ gives
		\[
		\pi(T;W,b)=\frac{T}{\varphi(W)\log T}+e_b(T),\qquad
		e_b(T)=o_W(T/\log T).
		\]
		The finite sum
		\[
		\rho_W(T)=\frac{\log T}{T}
		\sum_{\substack{b\bmod W\\(b,W)=1}}|e_b(T)|\longrightarrow0
		\]
		is independent of $a$. Therefore, simultaneously for all $a\in\PP$,
		\begin{align*}
			E_k(T;a)
			&\le2\pi(z)+\sum_{b\in\mathcal R_a}\pi(T;W,b)\\
			&\le2\pi(z)+\left(\frac{h_a}{\varphi(W)}+\rho_W(T)\right)
			\frac T{\log T}\\
			&\le2\pi(z)+\bigl(2\eta_k(z)+\rho_W(T)\bigr)\frac T{\log T}.
		\end{align*}
		
		\smallskip\noindent
		\emph{Step 3.}
		The inequality $\log(1-v)\le-v$ and the divergence of the sum of
		reciprocals of primes give
		\[
		\log\eta_k(z)
		\le-\sum_{\substack{3\le\ell\le z\\\ell\in\PP,\ \ell\nmid k}}
		\frac1{\ell-1}
		\le-\sum_{\substack{3\le\ell\le z\\\ell\in\PP,\ \ell\nmid k}}
		\frac1\ell\longrightarrow-\infty.
		\]
		Only finitely many primes divide $2k$, so their omission does not
		affect divergence; hence $\eta_k(z)\to0$. For each fixed $z$ the
		preceding bound, uniformly in $a$, yields
		\[
		\limsup_{T\to\infty}\sup_{a\in\PP}
		\frac{E_k(T;a)\log T}{T}\le2\eta_k(z).
		\]
		Letting $z\to\infty$ proves the lemma. The order of the limits
		is essential: $W$ is fixed while $T\to\infty$, and the resulting
		bound is independent of the prime shift $a$.
		
	\end{proof}
	
	\section{Proofs}
	\subsection{Separate limit sets and interval counts}
	We impose linear relations that prescribe the largest prime factors,
	select appropriate real proportions, and then use~\eqref{eq:GT} to
	count lattice points at which the forms take simultaneous prime values.
	
	\begin{proof}[Proof of Theorem~\ref{th:dense} on $\B$]
		\emph{Step 1.}
		Fix $m\ge1$. For primes $p\ne q$, the relations
		\[
		P(2p)=p,\qquad w(p^2q)=pP(p+q)^2
		\]
		suggest imposing $p+q=2mb$. Set
		\[
		p=a,\qquad q=2mb-a,\qquad
		\Phi_m(a,b)=\left(
		\frac{\max\{a,b\}}{\max\{a,2mb-a\}},
		\frac{b^2}{a(2mb-a)}\right).
		\]
		Given \(u>0\), we impose \(p_0/q_0=u\). Normalizing \(b_0=1\) and solving \(p_0+q_0=2mb_0\) yields the real centre
		\[
		(a_0,b_0)=\left(\frac{2mu}{u+1},1\right).
		\]
		At this centre, $p_0=2mu/(u+1)$ and
		$q_0=2m/(u+1)$ are positive, with
		\[
		\frac{p_0}{q_0}=u,
		\qquad
		\frac{b_0}{q_0}=\frac{u+1}{2m}.
		\]
		Substituting these ratios into the two
		coordinates of $\Phi_m$ gives
		\[
		\Phi_m(a_0,b_0)=F_m(u).
		\]
		
		Fix $\varepsilon>0$. Since the three forms
		$a$, $b$, and $2mb-a$ are strictly positive
		at $(a_0,b_0)$, and $\Phi_m$ is continuous
		in a neighbourhood of this point, we may choose a sufficiently small closed box \(K\) centred at \((a_0,b_0)\), with nonempty interior, on which all three forms remain positive and such that
		\[
		\|\Phi_m(a,b)-F_m(u)\|<\varepsilon
		\qquad ((a,b)\in K).
		\]
		By compactness and positivity, there exist
		constants $c,C>0$ such that
		\[
		c\le a,b,2mb-a\le C
		\qquad ((a,b)\in K).
		\]
		The norm is Euclidean, and all parameters
		and bounds are fixed independently of $N$.
		Since $\Phi_m$ is homogeneous of degree
		zero, while the three linear forms are
		homogeneous of degree one, we have
		\[
		\Phi_m(Ns,Nt)=\Phi_m(s,t)
		\qquad ((s,t)\in K).
		\]
		Consequently, for every $(a,b)\in NK$,
		writing $(a,b)=N(s,t)$ with $(s,t)\in K$
		gives
		\[
		\begin{aligned}
			cN &\le a,b,2mb-a\le CN,\qquad \|\Phi_m(a,b)-F_m(u)\|&<\varepsilon.
		\end{aligned}
		\]

		\smallskip\noindent
		\emph{Step 2.}
		Lemma~\ref{lem:three} and~\eqref{eq:GT} give
		\[
		\#\{(a,b)\in NK\cap\mathbb Z^2:a,b,2mb-a\in\PP\}
		\gg_{m,K}\frac{N^2}{(\log N)^3}.
		\]
		The degenerate case $p=q$ lies on a single line:
		\[
		p=q\ \Longleftrightarrow\ a=mb,\qquad
		\#\{(a,b)\in NK\cap\mathbb Z^2:a=mb\}
		=O_{m,K}(N)=o\!\left(\frac{N^2}{(\log N)^3}\right).
		\]
		After removing this line, the lower bound remains valid. 
		For every retained point, let $n=p^2q$.
		If $N>P(2m)/c$, then
		$b\ge cN>P(2m)$, and hence
		\[
		P(p+q)=P(2mb)=b.
		\]
		Since $P(2p)=p$, it follows that
		\begin{equation}\label{eq:Bidentities}
			w(n)=pb^2,
			\qquad
			(R(n),S(n))=\Phi_m(a,b).
		\end{equation}
		Thus $n\in\B$ and
		$\|(R(n),S(n))-F_m(u)\|<\varepsilon$.
		The map $(a,b)\mapsto p^2q$ is injective on the retained points:
		unique factorization recovers $p=a$ and $q$, and then
		$b=(p+q)/(2m)$. 
		Moreover,
		\[
		c^3N^3\le n=p^2q\le C^3N^3.
		\]
		In particular, all constructed integers
		tend to infinity with $N$.
		
		For sufficiently large $X$, take
		$N=X^{1/3}/C$. Then every constructed
		integer satisfies $n\le X$.
		By the injectivity above and the
		prime-point lower bound, we obtain
		\[
		\#\{n\le X:n\in\B,\ \|(R(n),S(n))-F_m(u)\|<\varepsilon\}
		\gg_{m,K}\frac{N^2}{(\log N)^3}
		\asymp_C\frac{X^{2/3}}{(\log X)^3},
		\]
		where $\log N=\frac13\log X+O_C(1)$.
		
		\smallskip\noindent
		\emph{Step 3.}
		For $\alpha\in(0,1)$, choose $m$ so that
		$\alpha>1/(2m-1)$ and take $u=\alpha$.
		Then $(u+1)/(2m)<u<1$, hence
		$r_m(u)=\alpha$. Given a nonempty open interval
		$I\subset(0,1)$, choose $\alpha\in I$
		and $\varepsilon>0$ such that
		$(\alpha-\varepsilon,\alpha+\varepsilon)
		\subset I$.
		Since
		\[
		|R(n)-\alpha|
		\le \|(R(n),S(n))-F_m(u)\|,
		\]
		the construction in Step 2 gives the
		required lower bound for $R(n)\in I$.
		
		For $\beta>0$, choose $m$ so that
		$m^{-2}<\beta$. Since $g_m$ is continuous
		on $(0,\infty)$ and
		\[
		g_m(1)=m^{-2}<\beta,
		\qquad
		\lim_{u\to\infty}g_m(u)=\infty,
		\]
		the intermediate value theorem gives
		some $u>1$ with $g_m(u)=\beta$.
		Given a nonempty open interval
		$J\subset(0,\infty)$, choose $\beta\in J$
		and $\varepsilon>0$ such that
		$(\beta-\varepsilon,\beta+\varepsilon)
		\subset J$.
		Since
		\[
		|S(n)-\beta|
		\le \|(R(n),S(n))-F_m(u)\|,
		\]
		Step 2 yields the required lower bound
		for $S(n)\in J$.

	\end{proof}
	
	\begin{proof}[Proof of Theorem~\ref{th:dense} on $\C$]
		\smallskip\noindent
		\emph{Step 1.}
		Fix $m\ge1$ and impose
		\[
		p+q=2mx,\qquad p+r=2my,\qquad q+r=2z.
		\]
		Solving for $p,q,r$ yields
		\[
		p=mx+my-z,\qquad q=mx-my+z,\qquad r=-mx+my+z.
		\]
		The six forms $x,y,z,p,q,r$ are precisely those in
		Lemma~\ref{lem:six}.
		Thus, requiring $x,y,z,p,q,r$ to be
		prime amounts to requiring simultaneous
		prime values of the six linear forms
		in Lemma~\ref{lem:six}.
		Suppose these six forms take prime
		values with $x,y,z>P(2m)$. Then
		\[
		\begin{aligned}
			P(p+q)&=P(2mx)=x,\\
			P(p+r)&=P(2my)=y,\\
			P(q+r)&=P(2z)=z.
		\end{aligned}
		\]
		Consequently,
		\begin{equation}\label{eq:Cidentities}
			w(pqr)=xyz.
		\end{equation}
		On the region where all six forms
		$x,y,z,p,q,r$ are positive, define
		the real quotient map
		\[
		\Psi_m(x,y,z)=
		\left(
		\frac{\max\{x,y,z\}}{\max\{p,q,r\}},
		\frac{xyz}{pqr}
		\right),
		\]
		where $p,q,r$ denote the linear
		forms above. This map is continuous
		and homogeneous of degree zero:
		\[
		\Psi_m(Nx,Ny,Nz)=\Psi_m(x,y,z)
		\qquad (N>0).
		\]
		In particular, whenever the six forms
		take prime values as above and
		$p,q,r$ are pairwise distinct, identity \eqref{eq:Cidentities} gives
		\[
		(R(pqr),S(pqr))=\Psi_m(x,y,z).
		\]
		
		\smallskip\noindent
		\emph{Step 2.}
		For fixed $m\ge1$ and $t>1/(2m)$,
		consider the real centre
		\[
		(x,y,z)=(t,t,1).
		\]
		At this point,
		\[
		(p,q,r)=(2mt-1,1,1),
		\]
		so all six forms are strictly positive.
		Moreover,
		\[
		\Psi_m(t,t,1)=
		\left(
		\frac{\max\{t,1\}}{\max\{2mt-1,1\}},
		\frac{t^2}{2mt-1}
		\right).
		\]
		
		To prescribe the first coordinate,
		let $I\subset(0,1)$ be a nonempty
		open interval and choose $\alpha\in I$.
		Choose $m$ sufficiently large that
		\[
		\alpha>\frac{1}{2m-1},
		\qquad
		t=\frac{1+\alpha^{-1}}{2m}.
		\]
		Then $1/m<t<1$ and
		$2mt-1=\alpha^{-1}>1$. Hence
		\[
		\frac{\max\{t,1\}}
		{\max\{2mt-1,1\}}
		=\frac{1}{2mt-1}
		=\alpha.
		\]
		To prescribe the second coordinate,
		let $J\subset(0,\infty)$ be a nonempty
		open interval and choose $\beta\in J$.
		Choose $m$ so that $m^{-2}<\beta$,
		and put
		\[
		h_m(t)=\frac{t^2}{2mt-1},
		\qquad t>\frac{1}{2m}.
		\]
		Since $h_m$ is continuous and
		\[
		h_m(1/m)=m^{-2}<\beta,
		\qquad
		\lim_{t\downarrow 1/(2m)}h_m(t)=\infty,
		\]
		the intermediate value theorem gives
		some $t\in(1/(2m),1/m)$ such that
		$h_m(t)=\beta$. In either case, fix the resulting
		parameters $m$ and $t$. Since all
		six forms are strictly positive at
		$(t,t,1)$ and $\Psi_m$ is continuous
		near this point, we may choose a
		sufficiently small closed box $K$
		centred at $(t,t,1)$, with nonempty
		interior, such that all six forms
		remain positive and the selected
		coordinate of $\Psi_m$ lies in $I$
		or $J$, respectively, throughout $K$.
		
		By compactness, there exist constants
		$c,C>0$ such that
		\[
		c\le x,y,z,p,q,r\le C
		\qquad ((x,y,z)\in K).
		\]
		All parameters and bounds are fixed
		independently of $N$. By homogeneity,
		the selected coordinate condition
		continues to hold on $NK$, and
		\[
		cN\le x,y,z,p,q,r\le CN
		\qquad ((x,y,z)\in NK).
		\]

		\smallskip\noindent
		\emph{Step 3.}
		Although $x=y$ and $q=r$ at the real centre, the box $K$ has
		nonempty interior in $\mathbb R^3$. To ensure that both $pqr$ and $w(pqr)$ belong to $\C$, we exclude
		equalities among the input primes
		$p,q,r$ and among the output primes $x,y,z$. Since
		\[
		p-q=2(my-z),\qquad p-r=2(mx-z),\qquad q-r=2m(x-y),
		\]
		the excluded set is the union $H$ of the hyperplanes
		\[
		my=z,\qquad mx=z,\qquad x=y,\qquad x=z,\qquad y=z.
		\]
		These are proper hyperplanes (some may coincide when $m=1$).
		By the hyperplane estimate established above,
		\[
		\#(H\cap NK\cap\mathbb Z^3)
		=O_{m,K}(N^2)
		=o\!\left(\frac{N^3}{(\log N)^6}\right).
		\]
		Combining Lemma~\ref{lem:six} with~\eqref{eq:GT} and removing $H$,
		we obtain
		\[
		\#\{(x,y,z)\in(NK\setminus H)\cap\mathbb Z^3:
		x,y,z,p,q,r\in\PP\}
		\gg_{m,K}\frac{N^3}{(\log N)^6}.
		\]
		For every retained point, Step 2 gives
		$cN\le x,y,z,p,q,r\le CN$, and hence
		\[
		c^3N^3\le n=pqr\le C^3N^3.
		\]
		If $N>P(2m)/c$, then $x,y,z>P(2m)$. Since the excluded
		hyperplanes ensure that $p,q,r$ and $x,y,z$ are each pairwise
		distinct, identity~\eqref{eq:Cidentities} gives
		\[
		n\in\C,\qquad w(n)=xyz\in\C,\qquad
		(R(n),S(n))=\Psi_m(x,y,z).
		\]
		For the selected interval, write $(Q,U)=(R,I)$ or $(S,J)$.
		By the choice of $K$ in Step 2 and the homogeneity of $\Psi_m$,
		every retained point satisfies $Q(n)\in U$.
		
		Each integer $n\in\C$ has at most $3!=6$ orderings of its prime
		factors. For each ordering, the identities
		\[
		x=\frac{p+q}{2m},\qquad
		y=\frac{p+r}{2m},\qquad
		z=\frac{q+r}{2}
		\]
		determine the prime point uniquely. Thus each $n$ arises from
		at most six retained points. For sufficiently large $X$, take $N=X^{1/3}/C$. Then every
		constructed integer satisfies $n\le X$. Since each integer arises
		from at most six retained points and
		$\log N=\frac13\log X+O_C(1)$, the prime-point estimate yields
		\[
		\#\{n\le X:n,w(n)\in\C,\ Q(n)\in U\}
		\gg_U\frac{X}{(\log X)^6}.
		\]
		The parameters $m,t,K,c,C$ depend only on the selected interval
		and are fixed before $X\to\infty$. This proves both bounds on $\C$.
		
		\smallskip\noindent
		\emph{Step 4.}
		Fix $\mathcal S\in\{\B,\C\}$ and
		$Q\in\{R,S\}$. The interval counts
		established above imply
		\[
		\#\{n\le X:n\in\mathcal S,\ Q(n)\in U\}
		\longrightarrow\infty
		\]
		for every nonempty open interval $U$
		contained in the interior of the asserted limit-point set of $Q$.
		
		For any interior target $\theta$,
		choose open intervals $U_j$ containing
		$\theta$, lying in that interior and in
		$(\theta-1/j,\theta+1/j)$.
		Starting from $n_0=0$, choose successively
		\[
		n_j\in\mathcal S,\qquad
		n_j>\max\{j,n_{j-1}\},\qquad
		Q(n_j)\in U_j.
		\]
		Such choices are possible because the
		interval counts are unbounded. Thus
		$n_j\to\infty$ and $Q(n_j)\to\theta$.
		The endpoints are obtained similarly:
		taking $U_j=(0,1/j)$ gives
		$Q(n_j)\to0$, while for $Q=R$,
		taking $U_j=(1-1/j,1)$ gives
		$R(n_j)\to1$. Finally, Chen, Shi and
		Wu~\cite[Corollary 5]{CSW} proved that
		\[
		\limsup_{\substack{n\to\infty\\n\in\A}}
		R(n)=1.
		\]
		Since $\mathcal S\subseteq\A$ and
		both quotients are positive, no other
		finite limit points can occur.
		This proves the asserted limit-point sets.
		
	\end{proof}
	
	\subsection{The joint limit set on \texorpdfstring{$\B$}{B3}}
	\begin{proof}[Proof of Theorem~\ref{th:joint}]
		Put $V(n)=(R(n),S(n))$. Let $\mathcal E$ denote
		the right-hand side of~\eqref{eq:joint}, and let
		$\LL_{\B}^{\mathrm{odd}}$ denote the set of finite
		limit points of $V(p^2q)$ as $p^2q\to\infty$,
		where $p$ and $q$ are distinct odd primes.
		It suffices to prove
		\[
		\mathcal E\subseteq\LL_{\B}^{\mathrm{odd}}
		\subseteq\LL_{\B}\subseteq\mathcal E.
		\]
		The middle inclusion follows from the definitions.
		All norms below are Euclidean norms.
		
		\smallskip\noindent
		\emph{Step 1.}
		Fix $m\ge1$, $u>0$ and $\varepsilon>0$.
		By the construction in the proof of
		Theorem~\ref{th:dense} on $\B$, leading
		to~\eqref{eq:Bidentities}, we may choose
		a closed box $K$ with nonempty interior
		and a constant $c_{m,K}>0$ such that,
		for all sufficiently large $N$,
		\[
		\#\{(a,b)\in NK\cap\mathbb Z^2:
		a,b,2mb-a\in\PP,\ a\ne mb\}
		\gg_{m,K}\frac{N^2}{(\log N)^3}.
		\]
		Here the excluded line $a=mb$, corresponding
		to $p=q$, contains only $O_{m,K}(N)$ lattice
		points in $NK$. This remains negligible
		when $u=1$, although the centre of $K$
		then lies on this line.
		
		For each retained point, put
		$p=a$, $q=2mb-a$ and $n=p^2q$.
		The same construction gives
		\[
		\|V(n)-F_m(u)\|<\varepsilon,\qquad
		\min\{p,q\}\ge c_{m,K}N,\qquad
		n\ge c_{m,K}^3N^3.
		\]
		Consequently, for every $M>0$, taking $N$
		sufficiently large yields an integer
		$n=p^2q>M$, with $p\ne q$ odd primes,
		such that
		\[
		\|V(n)-F_m(u)\|<\varepsilon.
		\]
		
		Starting with $n_0=0$, apply this assertion
		successively with $\varepsilon=1/j$ and
		$M=\max\{j,n_{j-1}\}$ to obtain integers
		$n_j$ with odd prime factors satisfying
		\[
		n_j>\max\{j,n_{j-1}\},\qquad
		\|V(n_j)-F_m(u)\|<1/j.
		\]
		Thus $n_j\to\infty$ and
		$V(n_j)\to F_m(u)$. Since $u>0$ and
		$m\ge1$ were arbitrary, we conclude that
		\[
		\bigcup_{m\ge1}\Gamma_m
		\subseteq\LL_{\B}^{\mathrm{odd}}.
		\]

		\smallskip\noindent
		\emph{Step 2.}
		The limit-point set $\LL_{\B}^{\mathrm{odd}}$ is closed,
		since
		\[
		\LL_{\B}^{\mathrm{odd}}
		=\bigcap_{M\ge1}
		\overline{\{V(p^2q):p,q\text{ distinct odd primes},
			\ p^2q>M\}}.
		\]
		Here the closure is taken in $\mathbb R^2$.
		
		For $0<A\le1$, we have
		\[
		F_m(A)=\left(
		\max\left\{A,\frac{A+1}{2m}\right\},
		\frac{(A+1)^2}{4m^2A}
		\right)\longrightarrow(A,0)
		\qquad(m\to\infty).
		\]
		For $A=0$, since $(m+1)/(2m^2)\le1/m$,
		\[
		F_m(1/m)=
		\left(\frac1m,\frac{(m+1)^2}{4m^3}\right)
		\longrightarrow(0,0).
		\]
		For $B>0$, put $u_m=(4Bm^2)^{-1}$.
		Then $u_m<1$ for all sufficiently large $m$, and
		\[
		F_m(u_m)=\left(
		\max\left\{\frac1{4Bm^2},
		\frac{1+1/(4Bm^2)}{2m}\right\},
		B(1+u_m)^2
		\right)\longrightarrow(0,B).
		\]
		Since $\Gamma_m\subseteq\LL_{\B}^{\mathrm{odd}}$
		for every $m\ge1$, closedness now gives
		\[
		([0,1]\times\{0\})\cup
		(\{0\}\times[0,\infty))
		\subseteq\LL_{\B}^{\mathrm{odd}}.
		\]
		Together with Step~1, this proves
		$\mathcal E\subseteq\LL_{\B}^{\mathrm{odd}}$.
		No uniformity in $m$ is needed, since Step~1
		applies to each fixed $m$ and closedness
		allows us to pass to limits as $m\to\infty$.
		
		For the reverse inclusion, we will also need
		the section with first coordinate $1$.
		For $u\ge1$ and $m\ge1$,
		\[
		u+1\le2u\le2mu,\qquad
		F_m(u)=(1,g_m(u)).
		\]
		Moreover,
		\[
		g_m(1)=\frac1{m^2},\qquad
		g_m(u)-\frac1{m^2}
		=\frac{(u-1)^2}{4m^2u}\ge0,
		\qquad
		\lim_{u\to\infty}g_m(u)=\infty.
		\]
		By continuity,
		$g_m([1,\infty))=[m^{-2},\infty)$.
		For every $B>0$, we may choose $m\ge1$
		with $m^{-2}\le B$, and then some $u\ge1$
		satisfies $g_m(u)=B$. Consequently,
		\begin{equation}\label{eq:vertical}
			\{1\}\times(0,\infty)
			\subseteq\bigcup_{m\ge1}\Gamma_m.
		\end{equation}
		The remaining point $(1,0)$ belongs to
		the horizontal axis.
		
		\smallskip\noindent
		\emph{Step 3.}
		Suppose $n_j=p_j^2q_j\to\infty$, where $p_j\ne q_j$
		are odd primes, and $V(n_j)\to(A,B)$.
		Then $A,B\ge0$. Put
		\[
		s_j=P(p_j+q_j),\qquad
		u_j=\frac{p_j}{q_j},\qquad
		m_j=\frac{p_j+q_j}{2s_j}.
		\]
		As shown before~\eqref{eq:curves}, we have
		$2s_j\mid p_j+q_j$, so $m_j$ is a positive integer.
		Using $P(2p_j)=p_j$ and~\eqref{eq:curves}, we obtain
		\[
		w(n_j)=p_js_j^2,\qquad
		\frac{s_j}{q_j}=\frac{u_j+1}{2m_j},\qquad
		V(n_j)=F_{m_j}(u_j).
		\]
		We may pass to subsequences throughout, since
		the limit $(A,B)$ remains unchanged.
		
		If $(m_j)$ is bounded, pass to a subsequence
		on which $m_j=m$ for some fixed $m\ge1$.
		Since $S(n_j)\to B<\infty$, there exists $C>0$
		such that
		\[
		\frac{u_j+2+u_j^{-1}}{4m^2}
		=S(n_j)\le C.
		\]
		Consequently,
		\[
		\frac1{4m^2C}\le u_j\le4m^2C.
		\]
		Passing to a further subsequence, we have
		$u_j\to u\in(0,\infty)$. By continuity of $F_m$,
		\[
		(A,B)=F_m(u)\in\Gamma_m\subseteq\mathcal E.
		\]
		
		If $(m_j)$ is unbounded, pass to a subsequence
		such that $m_j\to\infty$. Since
		$0<u_j/(1+u_j)<1$, we may pass further so that
		\[
		\frac{u_j}{1+u_j}\longrightarrow v\in[0,1].
		\]
		We distinguish three cases.
		
		\emph{(i) $0<v<1$.}
		Then $u_j\to u=v/(1-v)\in(0,\infty)$, so
		$u_j+u_j^{-1}=O(1)$. Hence
		\[
		S(n_j)=\frac{u_j+2+u_j^{-1}}{4m_j^2}
		\longrightarrow0,
		\]
		while
		\[
		R(n_j)
		=\frac{\max\{u_j,(u_j+1)/(2m_j)\}}
		{\max\{u_j,1\}}
		\longrightarrow\frac{u}{\max\{u,1\}}.
		\]
		Thus $(A,B)\in[0,1]\times\{0\}\subseteq\mathcal E$.
		
		\emph{(ii) $v=0$.}
		Then $u_j\to0$. Since $m_j\to\infty$,
		eventually $u_j<1$, and
		\[
		R(n_j)
		=\max\left\{u_j,\frac{u_j+1}{2m_j}\right\}
		\longrightarrow0.
		\]
		As $B$ is finite and nonnegative, we obtain
		\[
		(A,B)=(0,B)\in
		\{0\}\times[0,\infty)\subseteq\mathcal E.
		\]
		
		\emph{(iii) $v=1$.}
		Then $u_j\to\infty$. Eventually $u_j\ge1$,
		so $R(n_j)=1$ and hence $A=1$.
		If $B=0$, then
		$(A,B)=(1,0)\in[0,1]\times\{0\}$.
		If $B>0$, \eqref{eq:vertical} gives
		\[
		(A,B)=(1,B)\in\bigcup_{m\ge1}\Gamma_m.
		\]
		Thus $(A,B)\in\mathcal E$ in either case.
		
		We conclude that every finite limit point
		arising from distinct odd prime factors belongs
		to $\mathcal E$. Therefore
		\[
		\LL_{\B}^{\mathrm{odd}}\subseteq\mathcal E.
		\]

		\smallskip\noindent
		\emph{Step 4.}
		Suppose $n_j=p_j^2q_j\to\infty$, where $p_j\ne q_j$
		are primes, and $V(n_j)\to(A,B)$.
		Passing to an infinite subsequence, we may assume
		that either both primes are odd, $p_j=2$, or $q_j=2$
		throughout the subsequence. The first case was
		treated in Step~3.
		
		If $p_j=2$, then $n_j=4q_j\to\infty$, so
		$q_j\to\infty$. Put $s_j=P(q_j+2)\ge3$.
		Since $q_j$ is an odd prime, we have
		\[
		w(n_j)=2s_j^2,\qquad
		R(n_j)=\frac{s_j}{q_j},\qquad
		S(n_j)=\frac{s_j^2}{2q_j}.
		\]
		Consequently,
		\[
		R(n_j)^2=\frac{2S(n_j)}{q_j}
		\longrightarrow0,
		\]
		because $S(n_j)\to B<\infty$ and
		$q_j\to\infty$. Thus $A=0$, and hence
		\[
		(A,B)\in\{0\}\times[0,\infty)
		\subseteq\mathcal E.
		\]
		
		If $q_j=2$, then $n_j=2p_j^2\to\infty$,
		so $p_j\to\infty$. Put
		$s_j=P(p_j+2)\le p_j+2$.
		Using $P(2p_j)=p_j$, we obtain
		\[
		w(n_j)=p_js_j^2,\qquad
		1\le R(n_j)
		=\frac{\max\{p_j,s_j\}}{p_j}
		\le1+\frac2{p_j}\longrightarrow1.
		\]
		Thus $A=1$. If $B=0$, then
		$(A,B)=(1,0)\in[0,1]\times\{0\}$.
		If $B>0$, \eqref{eq:vertical} gives
		\[
		(A,B)=(1,B)\in\bigcup_{m\ge1}\Gamma_m.
		\]
		Hence $(A,B)\in\mathcal E$ in either case.
		
		Combining Steps~1--4, we have proved
		\[
		\mathcal E\subseteq\LL_{\B}^{\mathrm{odd}}
		\subseteq\LL_{\B}\subseteq\mathcal E.
		\qedhere
		\]

	\end{proof}
	
	\smallskip\noindent
	\emph{Remark.}
	Theorem~\ref{th:joint} gives an explicit description
	of the section of $\LL_{\B}$ at first coordinate $1/2$.
	By~\eqref{eq:joint}, the axis components contribute
	only $(1/2,0)$.
	
	For $u\ge1$, we have $r_m(u)=1$. For $0<u\le1$,
	\[
	r_m(u)=\max\left\{u,\frac{u+1}{2m}\right\}.
	\]
	If $m=1$, then $r_1(u)=(u+1)/2>1/2$.
	For $m\ge2$, the equality $r_m(u)=1/2$ holds
	if and only if $u=1/2$. In this case,
	\[
	g_m(1/2)
	=\frac{(1/2+1)^2}{4m^2(1/2)}
	=\frac9{8m^2}.
	\]
	Consequently,
	\[
	\LL_{\B}\cap(\{1/2\}\times[0,\infty))
	=\{(1/2,0)\}\cup
	\left\{\left(\frac12,\frac9{8m^2}\right):
	m\ge2\right\}.
	\]
	In particular, $(1/2,1)\notin\LL_{\B}$,
	although $1/2$ and $1$ are finite limit points
	of $R(n)$ and $S(n)$, respectively.

	\subsection{Relative density}
	
	\begin{proof}[Proof of Theorem~\ref{th:distribution}]
		Fix $\alpha>0$, and put
		\[
		E_{\mathcal S}(X)
		=|\{n\le X:n\in\mathcal S,\ R(n)>\alpha\}|
		\qquad(\mathcal S\in\{\A,\B,\C\}).
		\]
		We use the counting asymptotics of Chen, Shi
		and Wu~\cite[Lemma 6.1]{CSW}:
		\begin{equation}\label{eq:counts}
			\A(X)\sim\C(X)
			\sim\frac{X(\log\log X)^2}{2\log X},
			\qquad
			\B(X)\sim b_0\frac X{\log X},
			\qquad
			b_0=\sum_{p\in\PP}\frac1{p^2}.
		\end{equation}
		The constant $b_0$ is positive and finite.
		Choose a fixed constant $c>\max\{2,2/\alpha\}$,
		and put
		\[
		y=X^{1/10},\qquad
		U=X^{4/5}=\frac X{y^2}.
		\]
		Then $cy>y$, $\alpha cy>2y$, and
		$U/(cy)=X^{7/10}/c\to\infty$.
		
		\smallskip\noindent
		\emph{Step 1.}
		By the prime number theorem,
		\[
		\varepsilon(U):=
		\sup_{T\ge U}
		\left|\frac{\pi(T)\log T}{T}-1\right|
		\longrightarrow0
		\qquad(U\to\infty).
		\]
		Hence, uniformly for $T\ge U$,
		\[
		\pi(T)\ge
		(1-\varepsilon(U))\frac T{\log T}.
		\]
		The prime number theorem and partial summation
		also give (see~\cite[Chapters~2 and 3]{MV})
		\[
		\sum_{p\le v}\frac1p\sim\log\log v
		\qquad(v\to\infty).
		\]
		Define
		\[
		A_y=\sum_{p<q\le y}\frac1{pq},
		\qquad
		b_y=\sum_{p\le y}\frac1{p^2}.
		\]
		Then we get $b_y\longrightarrow b_0$ and
		\begin{align*}
			2A_y
			=\left(\sum_{p\le y}\frac1p\right)^2-b_y
			\sim(\log\log y)^2
			\sim(\log\log X)^2.
		\end{align*}
		Here
		$\log\log y=\log\log X+\log(1/10)$.
		Consequently, by~\eqref{eq:counts},
		\[
		\frac{XA_y}{\C(X)\log X}
		\longrightarrow1.
		\]

		\smallskip\noindent
		\emph{Step 2.}
		Put
		\[
		\mathcal T_C(X)
		=\{pqr\le X:p<q\le y,\ r>cy,\ p,q,r\in\PP\}.
		\]
		Since $p<q\le y<cy<r$, unique factorization
		ensures that each integer is counted only once.
		For $p<q\le y$, we have
		\[
		U\le\frac X{pq}\le X,
		\]
		and hence, for all sufficiently large $X$,
		\[
		\pi\left(\frac X{pq}\right)
		\ge(1-\varepsilon(U))
		\frac{X}{pq\log(X/(pq))}
		\ge(1-\varepsilon(U))\frac{X}{pq\log X}.
		\]
		Since $U>cy$ for sufficiently large $X$,
		\begin{align*}
			|\mathcal T_C(X)|
			&=\sum_{p<q\le y}
			\left\{\pi\left(\frac X{pq}\right)-\pi(cy)\right\}\\
			&\ge(1-\varepsilon(U))
			\frac X{\log X}A_y-cy^3,
		\end{align*}
		where
		\[
		\#\{p<q\le y\}\,\pi(cy)
		\le y^2\cdot cy=cy^3.
		\]
		Since $c$ is fixed, \eqref{eq:counts} gives
		\[
		\frac{cy^3}{\C(X)}
		\sim\frac{2c\log X}
		{X^{7/10}(\log\log X)^2}
		\longrightarrow0.
		\]
		Together with the estimates in Step~1, this yields
		\[
		1\ge\frac{|\mathcal T_C(X)|}{\C(X)}
		\ge(1-\varepsilon(U))
		\frac{XA_y}{\C(X)\log X}
		-\frac{cy^3}{\C(X)}
		\longrightarrow1.
		\]
		Consequently,
		\[
		\frac{|\mathcal T_C(X)|}{\C(X)}
		\longrightarrow1,
		\qquad
		\frac{|(\C\cap[1,X])\setminus\mathcal T_C(X)|}
		{\C(X)}
		\longrightarrow0.
		\]

		\smallskip\noindent
		\emph{Step 3.}
		For $n=pqr\in\mathcal T_C(X)$, we have
		\[
		P(n)=r,\qquad
		P(w(n))=\max\{P(p+q),P(p+r),P(q+r)\},
		\qquad
		P(p+q)\le2y<\alpha r.
		\]
		Thus $R(n)>\alpha$ implies that, for some
		$a\in\{p,q\}$,
		\[
		s=P(r+a)>\alpha r,\qquad
		r+a=ks,\qquad
		1\le k=\frac{r+a}{s}<\frac2\alpha,
		\]
		where $k$ is a positive integer and the last
		inequality follows from $a<r$.
		
		Put $H=\lceil2/\alpha\rceil$, and define
		\[
		\delta(V)=
		\sup_{\substack{T\ge V,\ a\in\PP\\1\le k\le H}}
		\frac{E_k(T;a)\log T}{T}
		\qquad(V\ge3).
		\]
		The supremum is finite since
		$E_k(T;a)\le\pi(T)\ll T/\log T$ for $T\ge3$.
		As $H$ is fixed, Lemma~\ref{lem:shift} gives
		\[
		\delta(V)\longrightarrow0
		\qquad(V\to\infty).
		\]
		For $p<q\le y$, we have $X/(pq)\ge U$ and
		$\log(X/(pq))\ge(4/5)\log X$. Hence
		\begin{align*}
			E_k\left(\frac X{pq};a\right)
			&\le\delta(U)\frac{X}{pq\log(X/(pq))}\\
			&\le\frac54\delta(U)\frac{X}{pq\log X}
			\qquad(a\in\PP,\ 1\le k\le H).
		\end{align*}
		
		Let $E_C^*(X)$ denote the number of integers
		$n\in\mathcal T_C(X)$ satisfying $R(n)>\alpha$.
		Dropping the restrictions $r>cy$ and
		$s=P(r+a)$, and applying the union bound, we obtain
		\begin{align*}
			E_C^*(X)
			&\le\sum_{p<q\le y}\sum_{k=1}^H
			\left\{
			E_k\left(\frac X{pq};p\right)
			+E_k\left(\frac X{pq};q\right)
			\right\}\\
			&\le\frac{5H}{2}\delta(U)\frac X{\log X}A_y.
		\end{align*}
		Possible overlaps only enlarge this upper bound. Since all exceptional integers outside
		$\mathcal T_C(X)$ lie in
		$(\C\cap[1,X])\setminus\mathcal T_C(X)$,
		we have
		\begin{align*}
			0\le\frac{E_{\C}(X)}{\C(X)}
			&\le\frac{E_C^*(X)}{\C(X)}
			+1-\frac{|\mathcal T_C(X)|}{\C(X)}\\
			&\le\frac{5H}{2}\delta(U)
			\frac{XA_y}{\C(X)\log X}
			+1-\frac{|\mathcal T_C(X)|}{\C(X)}
			\longrightarrow0.
		\end{align*}
		Thus $E_{\C}(X)=o(\C(X))$.

		\smallskip\noindent
		\emph{Step 4.}
		For this part, the fixed-shift Selberg-sieve bound of Chen and
		Shi~\cite[the sieve lemma on p.~2087]{CSparents} suffices:
		for fixed $k\ge1$ and $p\in\PP$,
		\[
		\#\{s\le T:s\in\PP,\ ks-p\in\PP\}
		\ll_{k,p}\frac{T}{(\log T)^2}.
		\]	
		Fix $p\in\PP$, and define
		\[
		Q_p=\max\{p,p/\alpha\},\qquad
		e_p(T)=\#\{q\le T:q\in\PP,\ q\ne p,\ R(p^2q)>\alpha\}.
		\]
		For $q>Q_p$, we have $q\ne p$ and $p<\alpha q$.
		Since $P(2p)=p$ and $P(p^2q)=q$, it follows that
		\[
		w(p^2q)=pP(p+q)^2,\qquad
		R(p^2q)=\frac{\max\{p,P(p+q)\}}{q}.
		\]
		Hence $R(p^2q)>\alpha$ implies
		\[
		s=P(p+q)>\alpha q,\qquad p+q=ks
		\]
		for some integer $k\ge1$, and therefore
		\[
		1\le k=\frac{p+q}{s}
		<\frac{p+q}{\alpha q}
		<\frac2\alpha.
		\]
		Put
		$
		H=\lceil2/\alpha\rceil.
		$
		If also $q\le T$ and $T\ge p$, then
		\[
		s\le p+q\le2T,\qquad q=ks-p.
		\]
		For each fixed $k$, the latter identity determines $q$
		uniquely from $s$. Hence
		\begin{align*}
			e_p(T)
			&\le \pi(Q_p)+\sum_{k=1}^{H}
			\#\{s\le2T:s\in\PP,\ ks-p\in\PP\}\\
			&\ll_{p,\alpha}1+\frac{T}{(\log T)^2}
			=o_{p,\alpha}\!\left(\frac{T}{\log T}\right).
		\end{align*}
		
		Now fix $Y\ge2$, and define
		\[
		\mathcal B_{\le Y}(X)
		=\{p^2q\le X:p,q\in\PP,\ p\le Y,\ p\ne q\},
		\qquad
		b_Y=\sum_{p\le Y}\frac1{p^2}.
		\]
		Unique factorization shows that the parametrization $p^2q$
		with $p\ne q$ is injective. For fixed $Y$, the prime number
		theorem gives
		\[
		|\mathcal B_{\le Y}(X)|
		=\sum_{p\le Y}\pi\!\left(\frac{X}{p^2}\right)+O_Y(1)
		\sim b_Y\frac{X}{\log X}.
		\]
		Since $Y$ is fixed, for each $p\le Y$ we have
		$X/p^2\to\infty$ and hence, for sufficiently large $X$,
		$X/p^2\ge p$. Therefore,
		\[
		\sum_{p\le Y}e_p\!\left(\frac{X}{p^2}\right)
		=o_{Y,\alpha}\!\left(\frac{X}{\log X}\right).
		\]
		All exceptional integers outside $\mathcal B_{\le Y}(X)$
		lie among the remaining
		$\B(X)-|\mathcal B_{\le Y}(X)|$ integers.
		Therefore,
		\[
		0\le E_{\B}(X)
		\le \B(X)-|\mathcal B_{\le Y}(X)|
		+\sum_{p\le Y}e_p\!\left(\frac{X}{p^2}\right).
		\]
		For each fixed $Y$, dividing by $\B(X)$ and letting
		$X\to\infty$, we obtain, by~\eqref{eq:counts} and
		the preceding estimates,
		\[
		0\le\limsup_{X\to\infty}
		\frac{E_{\B}(X)}{\B(X)}
		\le1-\frac{b_Y}{b_0}.
		\]
		Finally, $b_Y\to b_0$ as $Y\to\infty$.
		Since the preceding inequality holds for every
		fixed $Y$, it follows that
		$E_{\B}(X)=o(\B(X))$.

		\smallskip\noindent
		\emph{Step 5.}
		Since $\A=\B\mathbin{\dot\cup}\C$, we have
		\[
		E_{\A}(X)=E_{\B}(X)+E_{\C}(X),
		\qquad
		\A(X)=\B(X)+\C(X).
		\]
		Consequently,
		\[
		0\le\frac{E_{\A}(X)}{\A(X)}
		\le\frac{E_{\B}(X)}{\B(X)}
		+\frac{E_{\C}(X)}{\C(X)}
		\longrightarrow0.
		\]
		For $0<\alpha\le1$, it follows that
		\[
		D_{\A}(\alpha)
		=\lim_{X\to\infty}
		\left(1-\frac{E_{\A}(X)}{\A(X)}\right)
		=1.
		\qedhere
		\]
	\end{proof}

\end{document}